\documentclass[11pt]{article}
\usepackage{caption}
\usepackage{pdfsync}
\usepackage{amssymb}
\usepackage{color}
\usepackage{amsmath,amsthm}
\usepackage{graphicx}
\usepackage{empheq}
\usepackage{indentfirst}
\newtheorem{definition}{Definition}[section]

\newtheorem{remark}{Remark}[section]
\newtheorem{lem}{Lemma}[section]

\def\f{\frac}
\def\pa{\partial}

\def\e{\eqref}

\def\vep{\varepsilon}

\def\i1n{i=1,\cdots,n}
\def\j1n{j=1,\cdots,n}
\def\ij1n{i,j=1,\cdots,n}

\newcommand{\be}{\begin{equation}}
\newcommand{\ee}{\end{equation}}
\newcommand{\beq}{\begin{equation*}}
\newcommand{\eeq}{\end{equation*}}
 \numberwithin{equation}{section}

\newtheorem{thm}{Theorem}[section]
\DeclareMathOperator*{\esssup}{ess\,sup}
\title{Well-posedness of  the extrusion model described by coupled hyperbolic systems with a free boundary}
\author{
Peipei SHANG\thanks{Department of Mathematics, Tongji University, Shanghai 200092, China.
E-mail: {\tt peipeishang@hotmail.com}.}
\ Mamadou DIAGNE\thanks{Department of Mechanical and Aerospace Engineering,
University of California, San Diego,
La Jolla, CA, 92093, USA.
 E-mail: {\tt mdiagne@eng.ucsd.edu}}
\
\ and Zhiqiang WANG\thanks{School of Mathematical Sciences, Fudan
University, Shanghai 200433, China. E-mail: {\tt wzq@fudan.edu.cn}. } 
}
\date{}
\begin{document}

\maketitle

\begin{abstract}
In this paper, we consider the well-posedness of the Cauchy problem for
a physical model of the extrusion process,
which is described by  two systems of conservation laws with a free boundary.
By suitable change of coordinates and fixed point argument,
we obtain the existence,  uniqueness and regularity of the weak solution to this
Cauchy problem.
\\
\end{abstract}
{\bf 2010 MR Subject Classification:}\quad
         35L65, 
         35Q79,  
         35R37  

\section{Introduction}
Balance equations provide the foundation for much
physical-based modeling in fluid dynamics. They are also the
starting point for developing qualitative understanding of
phenomenological observations in fluid mechanics, heat transfer,
mass transfer, and reaction engineering.

The mathematical analysis of mobile interfaces in the context of
moving boundary problems has been an active subject in the last
decades and their mathematical understanding continues to be an
important interdisciplinary tool for the scientific applications.
Such kind of partial differential equation (PDE) model  arises in many
applications devoted to  modeling of biological systems and
reaction diffusion processes which involve stefan problems,
crystal growth processes.  One can mention applications which
concern swelling nanocapsules \cite{Bouchemal06}, lyophilization
\cite{Velardi08}, cooking processes \cite{Purlis10},
freeze drying process modeling \cite{Daraoui2010}, mixing
systems (model of torus reactor including a well-mixed zone and a
transport zone), diesel oxidation catalyst \cite{Petit10}.

Concerning biological systems,  \cite{CHOI11} proposes an analysis
of the global existence of solutions to a coupled
parabolic-hyperbolic system with moving boundary representing cell
mobility. A similar type of nonlinear moving-boundary problem,
consisting of a hyperbolic equation  and a parabolic equation for
modeling blood flow through viscoelastic arteries \cite{CAN11}
and tumor growth \cite{CHEN13},  is studied in terms of
well-posedness. Mathematical study of coupled partial differential
equations  through an internal moving interface  is also
proposed in \cite{Muntean09} for parabolic systems and in
\cite{Coutand11,Borsche10} for hyperbolic systems, where the PDEs are defined in a
time-varying spatial domain.

Generally speaking, the key resolution for such problems is based on a suitable
change of coordinates which
transforms the system with moving interface  into a system
defined on  a fixed domain. Then, the results of  many studies
dedicated to  systems of conservation laws may be useful to
establish the existence, uniqueness, regularity and continuous
dependence of solutions.  For the well-posedness problems, we
refer to the works \cite{Stefano05,BressanBook,LeFlochBook,LiuYang} (and the
references therein) in the content of weak solutions to systems
(including scalar case) of conservation laws, and
\cite{LiBook94,LiYuBook} in the content of classical solutions to
general quasi-linear hyperbolic systems.
Recall that there exists a classical approach developed for fixed
interfaces, consisting in augmenting the hyperbolic system of
conservation laws with color functions
\cite{Godlewski_ESAIM05,Godlewski04} for numerical analysis. For both cases with fixed
or moving interfaces, the lack of physical models which express
clearly the interface structure and the coupling conditions can be
considered as the real challenge from modeling point of view.
An interfacial model which describes precisely the information that
are exchanged at the coupling region should be defined with
respect to the real physical constraints.

In this paper we  consider the well-posedness of the Cauchy problem for
a physical model of the extrusion process. The process model is
composed by heat and mass transport equations which are defined in
complementary time-varying spatial domains. The domains are
coupled by a moving interface whose dynamics is governed by an
ordinary differential equation (ODE) expressing the conservation of mass
in an extruder.
More detailed description of the model is given in
Section \ref{description-model}.
We mention that  the first result concerning the
mathematical analysis  of the extrusion model as transport equations
coupled via complementary time varying domains is proposed in
\cite{DiagneJESA11}, where the well-posedness for the linearized
model of the extruder is obtained by using perturbation
theory on the linear operator.

Our proof of the well-posedness of the  Cauchy
problem for the extrusion model relies on a  change of
coordinates and a fixed point argument (see {\bf Theorem \ref{thm-well}}).
To tackle with the difficulty caused by the moving interface,
we make suitable change of coordinates  on the spatial variables
so that the moving interface problem is normalized to a standard
fixed domain problem. In order to deal with the nonlinearities,
we use Banach fixed point  theorem  based on the characteristic method,
which enables to compute the solution numerically.
The $H^2$-regularity of the solution is proved as well (see {\bf Theorem
\ref{thm-regu}}), which is useful when one considers the asymptotic
stabilization of the corresponding closed-loop system with
feedback controls
(see, e.g., \cite{CBdA} for hyperbolic systems with boundary
feedback laws).

We point out here that the analysis of the  Cauchy problem for the extrusion model
is fundamental to numerical simulation of the physical process. Moreover,  it is also
the first step for further study on control of the model.
It is of particular interest to consider some control problems,
including  controllability and stabilization of the filling ratio, the net flow rate, the position of the interface,
the moisture and the temperature of the extrusion process.
These problems will be studied in some forthcoming papers.

The organization of this paper is as follows: First in Section
\ref{description-model}, we give a description of the extrusion
process model  which is derived from conservation laws. As
preliminaries, we make domain normalization by change of
coordinates on space variables
in Section \ref{normalization-model}. The main results
 ({\bf Theorems \ref{thm-well}, \ref{thm-regu}, \ref{thm-mt}}) concerning
the well-posedness and regularity of the normalized system are
presented in Section \ref{main-results},
while their proofs are given in Section  \ref{sec-proof 1}-\ref{sec-proof 2} respectively.
Finally, in Section
\ref{conclusion}, we give the conclusion of this paper as well as
some perspectives.


\section{Description of the extrusion process model} \label{description-model}

\begin{figure}[h]
\begin{center}
\includegraphics[scale=0.4]{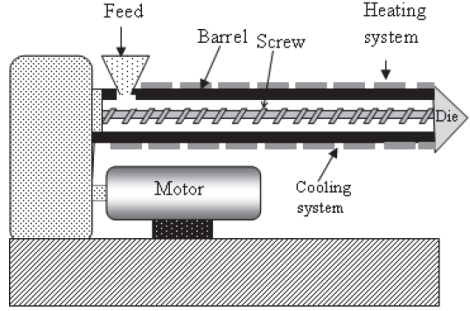}\vspace{0.005cm}
\caption{Schematic description of an extruder} \label{fig1}
\end{center}
\end{figure}

Extruders are designed to process highly viscous materials. They
are mainly used in the chemical industries for polymer processing
as well as in the food industries. An extruder is made of a
barrel, the temperature of which is regulated. One or two
Archimedean screws are rotating inside the barrel. The extruder is
equipped with a die where the material comes out of the process
(see Fig \ref{fig1}). The extruder is of particular interest
due to its modular geometry that allows the control of capacities
of the mixtures along the machine. Another interesting property of
the extruder is that the filling ratio along the axial direction
of the screws can be less than one in some part of the system
according to the screw configuration and the operating conditions.
For modeling purposes, the main phenomenon is obviously the fluid
flow which may be considered as highly viscous Newtonian or
non-Newtonian fluid flows interacting with heat transfer and
possibly chemical reactions. These processes occur within a
complex non-stationary volume delimited by the barrel and the
rotating screw. The most important part of the extruder is the screw
configuration which modulates extensively the mechanical energy.

\begin{figure}[!h]
\begin{center}
\includegraphics[scale=0.4]{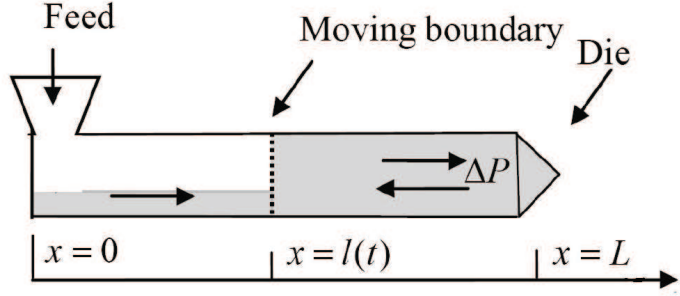}\vspace{0.005cm}
\caption{Bi-zone model of an extruder } \label{fig2}
\end{center}
\end{figure}

The complexity of screw geometric configuration in an extruder
make difficult the design a non-isothermal flow model
\cite{Booy80}. In \cite{Booy78}, an analysis of the flow in the
channel of co-rotating twin-screws at the same speed is developed
and the authors show how a reasonable flow analysis can be made by
writing a single screw extrusion process as an equivalent model.
In an extruder, the net flow at the die exit is mainly due to the
flow of the material in the longitudinal direction if one neglects the
clearance between the screw and the barrel and the vibrations
which can occur due to screws structure. Therefore, the flow
dynamics which is dominated by the convection effect in the
direction of screw axis is sufficient to represent the material
flow. This means that the transverse flow corresponding to a
recirculation of the material in the plane perpendicular to the
screw channel is neglected. From a macroscopic point of view, a
1D model describes clearly the material convection aspect in an
extruder. So, the material is driven from the feed to the die by
the pumping effect of the screw rotation. The geometric structure
of the die  influences the  transport along the extruder.
Therefore, the material is accumulated behind the die and fills
completely the available volume at this region. The spatial domain
where the extruder is completely filled is called the \emph{Fully
Filled Zone} ($FFZ$). The flow in this \emph{FFZ}
depends on the pumping capacity of the screw and also on the
pressure flow. The pressure gradient which appears due to the die
restriction is given by Navier-Stokes equation which provides a
mathematical model of the fluid motion. The extruder which is
initially empty, may also comprises a spatial domain that is not
completely filled by the material. This region  which corresponds
to a conveying region is called \emph{Partially Filled
Zone} ($PFZ$). In this domain, there is no pressure build-up. This
means that the pressure gradient is zero and the pressure is generally
equal to the air pressure inside the barrel. The transport
velocity of the material is controlled by the screw speed. These
two zones are coupled by an interface which characterizes the
spatial domains where the pressure gradient is null or not (and
accordingly where the extruder is partially or fully filled). The
mobile interface is assumed to be thin, i.e., reduced to a point. In
the sequel, the spatial domain of the extruder will be taken to be
the real interval $\left[0,\, L\right]$ where $L>0$ is the length
of the extruder. Let us denote by $l(t)\in\left[0,\, L\right]$ the
position of the thin interface, the domain of the \emph{PFZ}
is then $\left[0,\, l(t)\right]$ and the \emph{FFZ} is defined on $\left[l(t),\, L\right]$ ,
see Fig \ref{fig2}. The
interface is moving according to the volume of the material which is
accumulated in the \emph{FFZ}.  It is clear that the
interface which separates these two zones evolves as a function of
the difference between the feed and die rates. Finally, the extruder
model is composed by three interdependent dynamics which describe
the evolution of the material in the \emph{PFZ} and \emph{FFZ} and the evolution of the position of
the interface.

The extrusion model is based on this structural
decomposition of the extruder and  derived from the mass and the energy
balances  as in \cite{KULSH92,CHIN01}. The transport
equations describe the evolution of the filling ratio, the moisture
contains and the temperature for an extruded material. The problem of
coupled PDEs through an moving interface arises from the existence
of transport equations which are defined in two complementary
time-varying spatial domains denoted by $[0,\, l(t)]$
and  $[l(t),\, L]$.

\subsection{Physical definition of the parameters}
\begin{tabular}{ll|ll}
$L$       &\text{Extruder Length}    & $B$             &\text{Geometric parameter}              \\
$F_d$          &\text{Net forward mass flow rate}    & $K_d$    &\text{Geometric parameter} \\
$S_{ech}$  &\text{Exchange area between melt and barrel}   &  $V_{eff}$  &  \text{Effective volume}\\
$\alpha$ &\text{Heat exchange coefficient}   &  $S_{eff}$  &  \text{Effective area} \\ 
 {\color{red}$c_{o}$}      &\text{Specific heat capacity} & $\eta$  & \text{Melt viscosity} \\
 {\color{red}  $\beta_{o}$} & \text{Coefficient of viscous heat generation}     &  ${\color{red} \rho_{o} }$  &\text{Melt density}  \\
{\color{red} $\mu_p,\, \mu_f$} & \text{Viscous heat generation factor }  &   $\zeta$  &\text{Screw Pitch} \\
\end{tabular}

\subsection{The Partially Filled Zone ($PFZ$)  $[0,T]\times[0,l(t)]$}
For the $PFZ$, we consider the filling
ratio $f_{p}$, the moisture content $M_{p}$ and the temperature
$T_{p}$ as the state variables. The transport equations
associated with these variables are defined on $[0,T]\times[0,l(t)]$,
where $l(t)$ represents the moving
interface between the two zones. The filling ratio $f_p$ and the
moisture $M_p$  are strictly positive functions which are less than
one according to the modeling assumptions. The balance equations
express the convection through the rotation of the screw at
translational velocity $\alpha_p$, product on the pitch of the screw $\zeta$
and the rotation speed of the  screw $N(t)$. The source term
$\Omega_p$ which appears in the equation of $T_p$ groups
the heat produced by the rotational screw
(proportional to $N^{2}(t)$) and the heat exchange with the barrel
($T_b(t,x)$ is the distributed barrel temperature).
\begin{eqnarray}\label{orequ1}
\partial_t \begin{pmatrix}f_{p}\\
M_{p}\\
T_{p}\end{pmatrix}& = & -\alpha_p \partial_x\begin{pmatrix}f_{p}\\
M_{p}\\
T_{p}\end{pmatrix}+\begin{pmatrix}0\\
0\\
\Omega_p \end{pmatrix},\label{eqpfz}
\end{eqnarray}
where
\begin{align}
\alpha_p &=\zeta N(t) \\
\Omega_p & =  \frac{{\color{red} \mu_p \beta_{o}} \eta N^{2}(t)}{f_p (t,x){\color{red} \rho_{o} }V_{eff}{\color{red} c_{o}}}
+\frac{{\color{red} \zeta} S_{ech}\alpha}{{\color{red} \rho_{o} }V_{eff}{\color{red} c_{o}}}(T_{b}(t,x)-T_{p}(t,x)).
\label{eq:SourceBilanEnergie1}
\end{align}


\subsection{The Fully Filled Zone ($FFZ$) $[0,T]\times [l(t),L]$}

 For the $FFZ$, we consider  the moisture content $M_{f}$ and the
temperature $T_{f}$ as variables of the states. The transport equations
associated with these variables are defined on
$[0,T]\times [l(t),L]$, where $L$ is the length
of the extruder.  The transport velocity $\alpha_f$ is proportional to  the net flow rate $F_{d}(t)$ at the die.
Similarly as in the $PFZ$, the source term
$\Omega_f$  stands for  the heat produced by the rotational screw
and the heat exchange with the barrel in the $FFZ$.
\begin{eqnarray}\label{orequ2}
\partial_t\begin{pmatrix}M_{f}\\
T_{f}\end{pmatrix} & = & -\alpha_f \partial_x
\begin{pmatrix}M_{f}\\
T_{f}\end{pmatrix}+\begin{pmatrix}0\\
\Omega_f \end{pmatrix},\label{eqffz}
\end{eqnarray}
where
\begin{align}
\alpha_f  &=\frac{\zeta F_{d}(t)}{{\color{red} \rho_{o} }V_{eff}}
  \\
\Omega_f & =  \frac{{\color{red} \mu_f \beta_{o}}\eta N^{2}(t)}{{\color{red} \rho_{o} }V_{eff} {\color{red} {\color{red} c_{o}}}}
+\frac{{\color{red} \zeta}S_{ech}\alpha}{{\color{red} \rho_{o} }V_{eff} {\color{red} {\color{red} c_{o}}}}(T_{b}(t,x)-T_{f}(t,x)).
\label{eq:SourceBilanEnergie2}
\end{align}

\subsection{The interface}
Following \cite{KIM001,KIM002,KULSH92,CHIN01}, we assume that the two zones
are separated by an interface defined by the discontinuity of the
filling ratio. By assuming the continuity of pressure at the
interface $l(t)$, we  express the net flow rate $F_d$ as a function of
$l(t)$ and $N(t)$ as following:
\begin{align}
\label{ffd}
F_d(t)& =\frac{K_{d}}{\eta} \Delta P(t), \\
\Delta P(t)& :=P(t,L)-P_{0} =\frac{\eta {\color{red} \rho_{o} } V_{eff}N(t)(L-l(t))}{B{\color{red} \rho_{o} } +K_d (L-l(t))}.
\label{omega3}
  \end{align}
 The interface position is physically determined by the variation of pressure
 from the \emph{PFZ} to the \emph{FFZ} and  thus its dynamic is generated by the gradient
 of pressure which appears in the \emph{FFZ}. The equation  \eqref{omega3} is actually obtained by
 integrating the pressure gradient relation  from $l(t)$ to $L$ derived from momentum balance
in the \emph{FFZ} (where $B$ is a coefficient of pressure flow)
 \be
\partial_{x}P(t,x)=\eta\, \frac{{\color{red} \rho_{o} }V_{eff}N(t)-F_{d}(t)}{B{\color{red} \rho_{o} }}.
\label{pressure-flow}
\ee

We emphasis that without assuming a constant viscosity $\eta$,
which is  not distributed   with respect to the space variable,
along the extruder, the analytical solution of  \eqref{pressure-flow} can not be computed. 
Numerical simulations has been  performed by \cite{CHIN01}   considering the viscosity 
$\eta$ as a distributed function of moisture and temperature for steady state profile. 
In this case, the mass flow and the moisture are constant and the moving interface is stationary. 
The author uses iterative schemes assuming that  the pressure and the temperature profiles are 
known at the first step of the computation. The simulation shows the dynamics of the shaft power 
which is the power consumed by viscous dissipation in the $PFZ$ and the $FFZ$ and the power 
required to force the material through the die.

The interface dynamics which arises from a total mass balance  is
given by the following ODE
\be\label{eql}
\begin{cases}
\dot{l}(t)=F\left(l(t),N(t),f_p(t,l(t))\right),\\
l(0)=l^0,
\end{cases}
\ee where \be
F\left(l(t),N(t),f_p(t,l(t)\right)=\frac{F_{d}(t)-{\color{red} \rho_{o} }V_{eff}N(t)f_{p}(t,l(t))}{{\color{red} \rho_{o} }S_{eff}(1-f_{p}(t,l(t)))}.
\ee Recall that in the $PFZ$,  the filling ratio satisfies
 $f_{p}(t,x)<1,\:x\in(0,l(t))$ with $f_{p}(t,l(t))<1$ and in
 the $FFZ$, $f_{f}(t,x)=1,\:x\in (l(t),\,L)$.

 The moisture and the temperature are assumed to be continuous
at the interface
 \be\label{bound2}
\bigl(M_{p}(t,l(t)),\, T_{p}(t,l(t))\bigr)^{tr}
=\bigl(M_f(t,l(t)),\, T_f(t,l(t))\bigr)^{tr}.
\ee

\subsection{Initial and boundary conditions}
The initial conditions are given as the following
\begin{align}\label{initialc1}
(f_p(0,x),\, M_p(0,x),\, T_p(0,x))^{tr}&=(f^0_p(x),\, M^0_p(x),\, T^0_p(x))^{tr},\quad  x \in (0, l^0),\\
\label{initialc2}
(M_f(0,x),\, T_f(0,x))^{tr}&= (M^0_f(x),\, T^0_f(x))^{tr},\quad x\in (l^0,L).
\end{align}
%
The boundary conditions are given as the following
\be\label{bound1}
\bigl(f_{p}(t,0),\,
M_{p}(t,0),\,
T_{p}(t,0)\bigr)^{tr}=
\bigl(\frac{F_{in}(t)}{\theta(N(t))},\, M_{in}(t),\, T_{in}(t)\bigr)^{tr},
\ee
where
$F_{in}(t)$ denotes the feed rate and
\be \label{theta}
\theta(N(t))={\color{red} \rho_{o} }V_{eff}N(t).
\ee
Since the flow occurs in the direction of the screw channel,  $\theta$ is the maximum pumping capacity of the screw.

\section{Domain normalization of the extrusion process model}\label{normalization-model}

In this section, we will transform, by change of coordinates on space variables,
the original system with free boundary
(Cauchy problem  \eqref{orequ1}, \eqref{orequ2}, \eqref{eql}, \eqref{bound2}, \eqref{initialc1},
\eqref{initialc2} and \eqref{bound1}) to a normalized system with fixed boundary.

For the $PFZ$ zone, after change of variable
\beq
y=\f{x}{l(t)}
\eeq
 from $(0,l(t))$ onto the interval $(0,1)$ (see \cite{Diagne}),
we normalize system \eqref{eqpfz} to a new system defined on $Q:=(0,T)\times (0,1)$.
For the sake of simplicity, we still denote by $x$ the space variable instead of $y$,
the unknown functions by $(f_p,M_p,T_p)$, the velocity by $\alpha_p$
and the source term by $\Omega_p$. We have for all $(t,x)\in Q$ that
\begin{eqnarray}
\partial_t \begin{pmatrix}f_{p}(t,x)\\
M_{p}(t,x)\\
T_{p}(t,x)\end{pmatrix}+\alpha_p (x,N(t),l(t),f_p(t,1))\partial_x\begin{pmatrix}f_{p}(t,x)\\
M_{p}(t,x)\\
T_{p}(t,x)\end{pmatrix}=\begin{pmatrix}0\\
0\\
\Omega_p(t,x,T_p(t,x))   \end{pmatrix},   \label{eqpfz-normal}
\end{eqnarray}
with
\begin{align}\label{alphap-normal}
\alpha_p(x, N(t),l(t),f_p(t,1))& =\f{1}{l(t)}\bigl(\zeta N(t)-xF(l(t),N(t),f_p(t,1))\bigr),
  \\
\label{Omega_1}
\Omega_p(t,x,T_p(t,x)) &={\color{red} C_{o}}(T_p(t,x)-T_b(t,x))+g_p(t,x),
\end{align}
where
\be \label{g1}
{\color{red} C_{o}}=-\frac{{\color{red}\zeta}S_{ech}\alpha}{{\color{red} \rho_{o} }V_{eff}{\color{red} c_{o}}},
\quad
g_p(t,x)=\frac{  {\color{red} \mu_p \beta_{o}} \eta
N^{2}(t)}{f_p(t,x){\color{red} \rho_{o} }V_{eff}{\color{red} c_{o}}}.
\ee

For the $FFZ$ zone,
after change of variable
\beq
y=\f{x-l(t)}{L-l(t)}
\eeq
from $(l(t),L)$ onto the interval $(0,1)$ (see \cite{Diagne}).
system \eqref{eqffz} can be normalized to a new system defined on $Q$.
For the sake of simplicity, we still denote by $x$ the space variable instead of $y$,
the unknown functions by $(M_f,T_f)$, the velocity term by $\alpha_f$
and the source term by $\Omega_f$. Then we have for all $(t,x)\in Q$ that
\begin{eqnarray}
\partial_t\begin{pmatrix}M_{f}(t,x)\\
T_{f}(t,x)\end{pmatrix}+\alpha_f (x,N(t),l(t),f_p(t,1))\partial_x
\begin{pmatrix}M_{f}(t,x)\\
T_{f}(t,x)\end{pmatrix}=\begin{pmatrix}0\\
\Omega_f (t,x,T_f(t,x))\end{pmatrix},   \label{eqffz-normal}
\end{eqnarray}
with
\begin{align}\label{alphaf-normal}
\alpha_f (x,N(t),l(t),f_p(t,1))
&=\f{1}{L-l(t)}\Big(\frac{\zeta F_{d}(t)}{{\color{red} \rho_{o} }V_{eff}}+(x-1)F(l(t),N(t),f_p(t,1))\Big),
\\
\label{Omega_2}
\Omega_f (t,x,T_f(t,x))& = {\color{red} C_{o}} (T_f (t,x)-T_b(t,x))+g_f(t),
\end{align}
where
\begin{align}
F_d(t)  & =\frac{K_{d} {\color{red} \rho_{o} } V_{eff}N(t) (L-l(t))}{B{\color{red} \rho_{o} } +K_d (L-l(t))}, \\ \label{gf}
{\color{red} C_{o}} &  =-\frac{ {\color{red} \zeta}S_{ech}\alpha}{{\color{red} \rho_{o} }V_{eff} {\color{red} c_{o}}},\quad
g_f  (t)= \frac{  {\color{red} \mu_f \beta_{o}}  \eta
N^{2}(t)}{{\color{red} \rho_{o} }V_{eff} {\color{red} c_{o}}}.
\end{align}


The boundary conditions \eqref{bound2} and \eqref{bound1} can be rewritten as
\begin{align}\label{bound2-normal}
\bigl(M_f(t,0),\, T_f(t,0)\bigr)^{tr}&=\bigl(M_{p}(t,1),\,
T_{p}(t,1)\bigr)^{tr},\\
\label{bound1-normal}
\bigl(f_{p}(t,0),\,
M_{p}(t,0),\,
T_{p}(t,0)\bigr)^{tr}&=
\bigl(\frac{F_{in}(t)}{{\color{red} \rho_{o} }V_{eff}N(t)},\, M_{in}(t),\, T_{in}(t)\bigr)^{tr}.
\end{align}
In summary, we consider a coupled system
composed of an ODE for the moving interface
\be\label{eql-normal}
\begin{cases}
\dot{l}(t)=F\left(l(t),N(t),f_p(t,1)\right),\quad t\in(0,T),\\
l(0)=l^0,
\end{cases}
\ee
a transport equation for the filling ratio in the $PFZ$
\be\label{eq-filling-ratio-bound}
\begin{cases}
\partial_t f_p(t,x)+\alpha_p  (x,N(t),l(t),f_p(t,1))\partial_x f_p(t,x)=0,\quad (t,x)\in Q, \\
f_p(0,x)=f^{0}_p(x),\quad x\in(0,1),\\
f_p(t,0)=\displaystyle\frac{F_{in}(t)}{{\color{red} \rho_{o} }V_{eff}N(t)},\quad t\in(0,T),
\end{cases}
\ee
two transport equations for the moisture
\be\label{equation-moisture}
\begin{cases}
\pa_t M_p(t,x)+\alpha_p(x,N(t),l(t),f_p(t,1))\pa_x M_p(t,x)=0,\quad (t,x)\in Q,\\
\pa_t M_f(t,x)+\alpha_f(x,N(t),l(t),f_p(t,1))\pa_x M_f(t,x)=0,\quad (t,x)\in Q,\\
M_p(0,x)=M^0_p(x),\quad M_f(0,x)=M^0_f(x),\quad x\in(0,1),\\
M_p(t,0)=M_{in}(t),\quad M_f(t,0)=M_p(t,1),\quad t\in(0,T),
\end{cases}
\ee
and two transport equations for the temperature
\be\label{equation-temperature}
\begin{cases}
\pa_t T_p(t,x)+\alpha_p(x,N(t),l(t),f_p(t,1))\pa_x T_p(t,x)=\Omega_p(t,x,T_p(t,x)),\quad (t,x)\in Q,\\
\pa_t T_f(t,x)+\alpha_f(x,N(t),l(t),f_p(t,1))\pa_x T_f(t,x)=\Omega_f(t,x,T_f(t,x)),\quad (t,x)\in Q,\\
T_p(0,x)=T^0_p(x),\quad T_f(0,x)=T^0_f(x),\quad x\in(0,1),\\
T_p(t,0)=T_{in}(t),\quad T_f(t,0)=T_p(t,1),\quad t\in(0,T).
\end{cases}
\ee

\begin{remark}
Nonlinearity of the model is focused on the Cauchy problem \eqref{eql-normal}-\eqref{eq-filling-ratio-bound}
which is closed  for $(l,f_p)$.
With known values of $(l,f_p)$, Cauchy problem \eqref{equation-moisture}-\eqref{equation-temperature}
is linear with respect to the unknowns $(M_p,M_f,T_p,T_f)$.
\end{remark}

In the whole paper, unless otherwise specified,
we always assume that $l^0\in(0,L)$, $f^0_p\in W^{1,\infty}(0,1)$,
$M^0_p, T^0_p, M^0_f, T^0_f\in L^2(0,1)$,
$M_{in}, T_{in}\in L^2(0,T)$,
$F_{in}, N\in W^{1,\infty}(0,T)$ and $T_b\in L^2(Q)$.
For the sake of simplicity, we denote from now on
$\|f\|_{L^{\infty}}$ ($\|f\|_{W^{1,\infty}}$, $\|f\|_{L^2}$, resp.) as the $L^{\infty}$ ($W^{1,\infty}$, $L^2$, resp.)
norm of the function $f$ with respect to its variables.

\section{Main Results}\label{main-results}

In this section, we show the main results on the well-posedness of
the whole coupled system \eqref{eql-normal}-\eqref{equation-temperature}.
We first study the nonlinear Cauchy problem
\eqref{eql-normal}-\eqref{eq-filling-ratio-bound} since it is closed  for $(l,f_p)$,
then turn to linear Cauchy problem \eqref{equation-moisture}-\eqref{equation-temperature}
with known $(l,f_p)$.

Concerning the  Cauchy problem
\eqref{eql-normal}-\eqref{eq-filling-ratio-bound},
we have the following two theorems

\begin{thm}\label{thm-well}
Let $T>0$. Let $(l_e, N_e, f_{pe})$ be a constant equilibrium, i.e.,
\be\label{equil}
F(l_e,N_e,f_{pe})=0
\ee
with $0<f_{pe}<1$, $0<l_e<L$.
Assume that the compatibility condition at $(0,0)$ holds
\be \label{compati}
\f{F_{in}(0)}{{\color{red} \rho_{o} }V_{eff}N(0)}=f^0_p(0).
\ee
Then, there exists $\varepsilon_0$ (depending on $T$) such that
for any $\varepsilon\in(0,\varepsilon_0]$,
if
\begin{eqnarray}
\|f^{0}_p(\cdot)-f_{pe}\|_{W^{1,\infty}}
+ \|\frac{F_{in}(\cdot)}{{\color{red} \rho_{o} }V_{eff}N(\cdot)}- f_{pe}\|_{W^{1,\infty}}
+ \|N(\cdot)\!-\!N_e\|_{W^{1,\infty}}
+ |l^0\!-\!l_e|\!\leq \!\varepsilon, \label{G3}
\end{eqnarray}
Cauchy problem \eqref{eql-normal}-\eqref{eq-filling-ratio-bound}
admits a unique solution
$(l,f_p)\in W^{1,\infty}(0,T)\times W^{1,\infty}(Q)$,
and the following estimates hold
\begin{align} \label{est-fp}
\|f_p(\cdot,\cdot)-f_{pe}\|_{W^{1,\infty}}&\leq  C_{\varepsilon_0}\cdot\varepsilon,\\
\|l(\cdot)-l_e\|_{W^{1,\infty}}&\leq  C_{\varepsilon_0}\cdot\varepsilon,
\end{align}
where $C_{\varepsilon_0}$ is a constant depending on $\varepsilon_0$,
but independent of $\varepsilon$.
\end{thm}

\begin{thm}\label{thm-regu}
Under the assumptions of Theorem \ref{thm-well}, we assume furthermore that
$f^0_p(\cdot)\in H^2(0,1)$, $\displaystyle\f{F_{in}(\cdot)}{{\color{red} \rho_{o} }V_{eff}N(\cdot)}\in H^2(0,T)$,
and the compatibility condition at $(0,0)$ holds
\be \label{compati2}
(f^0_{p})_x(0)+\displaystyle\f{l(0)}{\zeta N(0)}\cdot\f{F'_{in}(0)N(0)-F_{in}(0) N'(0)}{{\color{red} \rho_{o} }V_{eff}N^2(0)}=0.
\ee
Then, there exists $\varepsilon_0$ (depending on $T$) such that for any $\varepsilon\in (0,\varepsilon_0]$, if
\begin{eqnarray}
\|f^{0}_p(\cdot)\!-\!f_{pe}\|_{H^2(0,1)}
  + \|\frac{F_{in}(\cdot)}{{\color{red} \rho_{o} }V_{eff}N(\cdot)} - f_{pe}\|_{H^2(0,T)}
  +\|N(\cdot)\!-\!N_e\|_{W^{1,\infty}} + |l^0\!-\!l_e|\!\leq \!\varepsilon,   \label{G333}
\end{eqnarray}
Cauchy problem \eqref{eql-normal}-\eqref{eq-filling-ratio-bound} has a unique solution
$(l,f_p)\in W^{1,\infty}(0,T)\times C^0([0,T];H^2(0,1))$
with the additional estimate
\be
\|f_p(\cdot,\cdot)-f_{p_e}\|_{C^0([0,T];H^2(0,1))}\leq  C_{\varepsilon_0}\cdot \varepsilon,
\ee
where $C_{\varepsilon_0}$ is a constant depending on $\varepsilon_0$,
but independent of $\varepsilon$.
\end{thm}

\begin{remark}
The solution in Theorem \ref{thm-well} or in Theorem \ref{thm-regu}
is often called semi-global solution since it exists on any preassigned time interval
$[0,T]$ if $(l,f_p)$ has some kind of smallness (depending on $T$), see \cite{Li2001, W2006}.
\end{remark}

\begin{remark}\label{rem}
We have the hidden regularity that $f_p \in C^0([0,1];H^2(0,T))$ in Theorem \ref{thm-regu}.
\end{remark}

For the  proof of Remark \ref{rem}, one can refer to \cite{CKW, SW}.

Concerning the moisture equation \eqref{equation-moisture}
and the temperature equation \eqref{equation-temperature},
we have the following theorem

\begin{thm}\label{thm-mt}
Under the assumptions of Theorem \ref{thm-well},
Cauchy problem \eqref{equation-moisture}-\eqref{equation-temperature}
admits a unique solution $(M_p,M_f,T_p,T_f)\in (C^0([0,T];L^2(0,1)))^4$,
and the following estimates hold
\begin{align} \label{est-mt}
\|M_p\|_{C^0([0,T];L^2(0,1))}&\leq   C\cdot\bigl(\|M^0_{p}\|_{L^2} +\|M_{in}\|_{L^2}\bigr),\\ \label{est-Tp}
\|T_p\|_{C^0([0,T];L^2(0,1))}&\leq  C\cdot\bigl(\|T^0_p\|_{L^2}+\|T_{in}\|_{L^2}+\|g_p\|_{L^2}\bigr),\\ \label{est-Mf}
\|M_f\|_{C^0([0,T];L^2(0,1))} &\leq  C\cdot \bigl(
 \|M^0_{p}\|_{L^2} +\|M_{in}\|_{L^2}+\|M^0_f\|_{L^2}\bigr),\\ \label{est-Tf}
\|T_f\|_{C^0([0,T];L^2(0,1))} &\leq  C\cdot \bigl(\|T^0_p\|_{L^2}+\|T_{in}\|_{L^2}+\|T^0_f\|_{L^2}+\|g_f\|_{L^2}\bigr),
\end{align}
where $g_p$ and $g_f$ are defined as \eqref{g1} and \eqref{gf} respectively and
$C>0$ is a constant.
\end{thm}


\section{Proof of Theorem \ref{thm-well}} \label{sec-proof 1}

In order to conclude Theorem \ref{thm-well}, it suffices to prove the following lemma on
local well-posedness of  Cauchy problem
\eqref{eql-normal}-\eqref{eq-filling-ratio-bound}.

\begin{lem}\label{lem-loc}
There exist $\varepsilon_1>0$ suitably small and $\delta=\delta(\varepsilon_1,\,
\|f^0_p(\cdot)-f_{pe}\|_{W^{1,\infty}},\, |l^0-l_e|) >0$, such that for any $\varepsilon\in(0,\varepsilon_1]$,
$f^0_p\in W^{1,\infty}(0,1)$, $F_{in}, N\in W^{1,\infty}(0,T)$,
$l^0\in (0,L)$ with
\be\label{initial-bound}
\|f^{0}_p(\cdot)\!-\!f_{pe}\|_{W^{1,\infty}}\!
   +\!\|\frac{F_{in}(\cdot)}{{\color{red} \rho_{o} }V_{eff}N(\cdot)}\!-\!f_{pe}\|_{W^{1,\infty}}\!
   +\!\|N(\cdot)\!-\!N_e\|_{W^{1,\infty}}\!
   +\!|l^0\!-\!l_e|\!\leq \!\varepsilon,
\ee
Cauchy problem \eqref{eql-normal}-\eqref{eq-filling-ratio-bound}
admits a unique local solution on $[0,\delta]$, which satisfies
the following estimates
\begin{align}\label{estimatefp}
\|f_p(t,\cdot)-f_{pe}\|_{W^{1,\infty}}
&\leq  C_{\varepsilon_1}\cdot\varepsilon,\quad \forall t\in[0,\delta],\\
\label{estimatel}
|l(t)-l_e|&\leq  C_{\varepsilon_1}\cdot\varepsilon,\quad \forall t\in [0,\delta],
\end{align}
where $C_{\varepsilon_1}$ is a constant depending on $\varepsilon_1$, but independent of $\varepsilon$.
\end{lem}

Let us first show how to conclude Theorem \ref{thm-well} from Lemma \ref{lem-loc}.
By Lemma \ref{lem-loc}, we take $\varepsilon_2 \in (0, \varepsilon_1]$
such that $C_{\varepsilon_1}\cdot\varepsilon_2\leq  \varepsilon_1$.
Then  for any $\varepsilon\in(0,\varepsilon_2]$ and any initial-boundary data such that
\eqref{initial-bound} holds,
Cauchy problem \eqref{eql-normal}-\eqref{eq-filling-ratio-bound} admits a unique solution on $[0,\delta]$.
Furthermore, one has
\begin{align}\label{n-use}
\|f_p(\delta,\cdot)-f_{pe}\|_{W^{1,\infty}}&\leq  C_{\varepsilon_1}\cdot\varepsilon\leq  \varepsilon_1,\\
\label{n-use2}
|l(\delta)-l_e|&\leq  C_{\varepsilon_1}\cdot \varepsilon\leq  \varepsilon_1.
\end{align}
By taking $(l(\delta),f_p(\delta,\cdot))$ as new initial data  and applying Lemma \ref{lem-loc} on $[\delta,2\delta]$,
the solution of Cauchy problem \eqref{eql-normal}-\eqref{eq-filling-ratio-bound} is extended to $[0,2\delta]$.
For fixed $T>0$, we can extend the local solution to Cauchy problem
 \eqref{eql-normal}-\eqref{eq-filling-ratio-bound} to $[0,T]$ eventually
by reducing the value of $\varepsilon_0$ and applying Lemma \ref{lem-loc}
in finite times (at most $[T/\delta]+1$ times).
Therefore, to conclude Theorem \ref{thm-well}, it remains to prove Lemma  \ref{lem-loc}.  \qed


\noindent
{\bf Proof of Lemma \ref{lem-loc}:}
The proof is divided into 4 steps.

\noindent \textbf{Step 1. Existence and uniqueness of $(l(\cdot),f_p(\cdot,1))$ by fixed point argument.}

Let $\varepsilon_1>0$ be such that
\be
0<\varepsilon_1<\min\{l_e,L-l_e,f_{pe},1-f_{pe}\}.
\ee
Denote
\begin{align}
\label{FW}
\|F\|_{W^{1,\infty}}:&=\sum_{|\alpha|\leq  1}\esssup_{\substack{0<x_1<L\\
N_e-\varepsilon_1<x_2<N_e+\varepsilon_1\\
0<x_3<1}}|D^{\alpha}F(x_1,x_2,x_3)|,\\
\Psi(t):&=(l(t),f_p(t,1)),\quad t\in[0,T].
\end{align}

\begin{figure}[htbp!]
\includegraphics[width=0.33\textwidth]{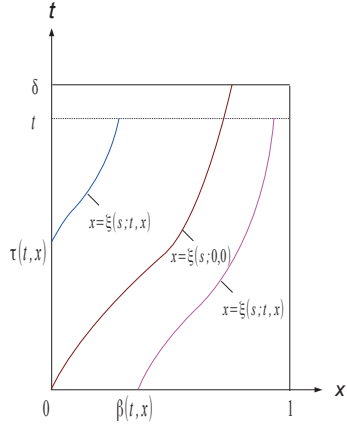}
\centering
\caption{The characteristics $\xi(s;t,x)$ \label{Fig3} }\par\vspace{0pt}
\end{figure}

For any given $\delta>0$ small enough (to be chosen later), we define
a domain candidate as a closed subset of $C^0([0,\delta])$ with respect to $C^0$ norm:
\begin{align}
\Omega_{\delta,\varepsilon_1}:=
\Big\{ \Psi \in C^0 ([0 ,\delta]):\, \Psi(0)=(l^0,f^0_{p}(1)),\
\|\Psi(\cdot)-(l_e,f_{pe})\|_{C^0([0,\delta])} \leq \vep_1 \Big\}.
\end{align}
We denote by $\xi(s;t,x)$, with $(s, \xi(s;t,x)) \in [0,t]\times [0,1]$
the characteristic curve passing through the point $(t,x) \in [0,\delta] \times [0,1]$ (see Fig \ref{Fig3}), i.e.,
\be\label{odexi2}
\begin{cases}
\displaystyle\frac{d \xi(s;t,x)} {ds}=\alpha_p\bigl(\xi(s;t,x),N(s),l(s),f_p(s,1)\bigr),\\
\xi(t;t,x)=x.
\end{cases}
\ee

Let us define a map $\mathfrak{F}:=(\mathfrak{F}_1,\mathfrak{F}_2)$, where
$\mathfrak{F} : \Omega_{\delta,\vep_1} \rightarrow  C^0 ([0 ,\delta])$, $\Psi\mapsto \mathfrak{F}(\Psi)$
as
\begin{align}\label{map1}
\mathfrak{F}_1(\Psi)(t)&:=l^0+\int_0^t F\bigl(l(s),N(s),f_p(s,1)\bigr)ds,\\
\label{map2}
\mathfrak{F}_2(\Psi)(t)&:=f^0_p(\xi(0;t,1)).
\end{align}
Solving the linear ODE \eqref{odexi2} with $\alpha_p$ given by \eqref{alphap-normal},
one easily gets for all $\delta $ small and all $0\leq  s\leq  t \leq \delta$ that
\be\label{solodexi2}
\xi(s;t,1)=e^{\int_s^t \f{F(l(\sigma),N(\sigma),f_p(\sigma,1))}{l(\sigma)}\, d\sigma}-\int_s^t
\f{\zeta N(\sigma)}{l(\sigma)}\cdot e^{\int^{\sigma}_s \f{F(l(s),N(s),f_p(s,1))}{l(s)}\, ds}\,d\sigma.
\ee
It is obvious that $\mathfrak{F}$ maps into $\Omega_{\delta,\vep_1}$ itself if
\be \label{delta01}
0<\delta< \min \Big\{T_0,T,\f{l_e-\vep_1}{\|F\|_{W^{1,\infty}}},\f{L-l_e-\vep_1}{\|F\|_{W^{1,\infty}}}\Big \},
\ee
where $T_0$ denotes the time when the characteristic curve $\xi(s;0,0)$ arrives at $x=1$, i.e., $\xi(T_0;0,0)=1$.

Now we prove that, if $\delta$ is small enough, $\mathfrak{F}$ is a contraction mapping on
$\Omega_{\delta,\vep_1}$ with respect to the $C^0$ norm.
Let $\Psi=(l,f_p),\bar \Psi=(\bar l,\bar f_p)\in \Omega_{\delta,\vep_1}$.
We denote by $\bar \xi(s;t,x)$ the corresponding characteristic curve  passing through $(t,x)$:
   \be
   \begin{cases}
\displaystyle\frac{d \bar \xi(s;t,x)} {ds}=\alpha_p\bigl(\bar \xi(s;t,x),N(s),\bar l(s), \bar f_p(s,1)\bigr),\\
\bar\xi(t;t,x)=x.
\end{cases}
   \ee
Similarly as \e{solodexi2},  one has for all $\delta $ small and all $0\leq  s\leq  t \leq \delta$ that
\be\label{solodebarxi2}
\bar\xi(s;t,1)=e^{\int_s^t \f{F(\bar l(\sigma),N(\sigma),\bar f_p(\sigma,1))}{\bar l(\sigma)}\, d\sigma}-\int_s^t
\f{\zeta N(\sigma)}{\bar l(\sigma)}\cdot e^{ \int^{\sigma}_s \f{F(\bar l(s),N(s),\bar f_p(s,1))}{\bar l(s)}\,
ds}\,d\sigma.
\ee

Therefore it holds for all $t\in[0,\delta]$ that
 \begin{align}
|\mathfrak{F}_1(\bar \Psi)(t)-\mathfrak{F}_1(\Psi)(t)|=&
\Big|\int_0^t F\bigl(\bar l(s),N(s),\bar f_p(s,1))\, ds-\int_0^t F\bigl(l(s),N(s),f_p(s,1))\, ds \Big |\nonumber\\
\leq  &\delta \|F\|_{W^{1,\infty}} \|\bar \Psi-\Psi\|_{C^0([0,\delta])}.   \label{difference1}
   \end{align}
On the other hand, it follows from \eqref{map2}, \eqref{solodexi2} and \eqref{solodebarxi2},   that for  all $t\in [0,\delta]$,
\begin{align}
      &|\mathfrak{F}_2(\bar \Psi)(t)-\mathfrak{F}_2(\Psi)(t)|\nonumber
        \\
    =&|f^0_p(\bar\xi(0;t,1))-f^0_p(\xi(0;t,1))|\nonumber
       \\
    \leq  & \|f^0_{px}\|_{L^{\infty}}\Big\{\Big|e^{\int_0^t \f{F(\bar l(\sigma),N(\sigma),\bar f_p(\sigma,1))}{\bar l(\sigma)}\, d\sigma}
                    -e^{\int_0^t \f{F(l(\sigma),N(\sigma),f_p(\sigma,1))}{l(\sigma)}\, d\sigma}\Big|\nonumber
       \\
      & \qquad \qquad + \int_0^t\Big|\f{\zeta N(\sigma)}{\bar l(\sigma)}\cdot
          e^{\int^{\sigma}_0 \f{F(\bar l(s),N(s),\bar f_p(s,1))}{\bar l(s)}\, ds}-\f{\zeta N(\sigma)}{l(\sigma)}\cdot
          e^{\int^{\sigma}_0 \f{F(l(s),N(s),f_p(s,1))}{l(s)}\, ds}\Big|\, d\sigma \Big\}.\nonumber
\end{align}
By \e{FW} and the fact that $\Psi,\bar \Psi \in \Omega_{\delta,\vep_1}$ and
$l(t)\geq l_e-\vep_1>0, |N(t)|\leq N_e+\vep_1 $, it follows  that for  all $t\in [0,\delta]$,
\begin{align}\label{difference22}
      &|\mathfrak{F}_2(\bar \Psi)(t)-\mathfrak{F}_2(\Psi)(t)|\nonumber
        \\
        \leq  & \|f^0_{px}\|_{L^{\infty}}\Big\{e^{\f{\delta\|F\|_{W^{1,\infty}}}{l_e-\vep_1}}
                        \int_0^t\Big|\f{F(\bar l(\sigma),N(\sigma),\bar f_p(\sigma,1))}{\bar l(\sigma)}
                             - \f{F(l(\sigma),N(\sigma),f_p(\sigma,1))}{l(\sigma)}\Big|\, d\sigma \nonumber
        \\
    &+   \f{\zeta (N_e+\vep_1)}{l_e-\vep_1} e^{\f{\delta\|F\|_{W^{1,\infty}}}{l_e-\vep_1}}
              \int_0^t \int_0^{\sigma}\Big|\f{F(\bar l(s),N(s),\bar f_p(s,1))}{\bar l(s)}- \f{F(l(s),N(s),f_p(s,1))}{l(s)}\Big|\, ds d\sigma
          \nonumber  \\
     &+  \zeta (N_e+\vep_1) e^{\f{\delta\|F\|_{W^{1,\infty}}}{l_e-\vep_1}}
                \int_0^t \Big|\f{1}{\bar l(\sigma)}- \f{1}{ l(\sigma)}\Big|\, d\sigma   \Big \}
          \nonumber \\
   \leq  & \|f^0_{px}\|_{L^{\infty}} e^{\f{\delta\|F\|_{W^{1,\infty}}}{l_e-\vep_1}}
         \Big\{ \delta\|F\|_{W^{1,\infty}} \Big (1+\f{\delta \zeta (N_e+\vep_1)}{l_e-\vep_1} \Big )
            \Big[  \f{ \|\bar \Psi-\Psi \|_{C^0([0,\delta])}  }{l_e-\vep_1} + \f{ \|\bar l- l\|_{C^0([0,\delta])}  }{(l_e-\vep_1)^2}  \Big]
         \nonumber  \\
     &+  \delta \zeta (N_e+\vep_1)    \f{ \|\bar l- l\|_{C^0([0,\delta])}} {(l_e-\vep_1)^2}      \Big \}
         \nonumber  \\
   \leq  &  \delta  e^{\f{T\|F\|_{W^{1,\infty}}}{l_e-\vep_1}}
         \Big\{ \Big (1+\f{T \zeta (N_e+\vep_1)}{l_e-\vep_1} \Big )
            \Big[  \f{\|F\|_{W^{1,\infty}}  }{l_e-\vep_1} + \f{ \|F\|_{W^{1,\infty}} }  {(l_e-\vep_1)^2}  \Big]
   +    \f{\zeta (N_e+\vep_1)  } {(l_e-\vep_1)^2}      \Big \}   \nonumber
     \\
      & \qquad \qquad \cdot   \|f^0_p(\cdot)-f_{pe}\|_{W^{1,\infty}}  \|\bar \Psi-\Psi \|_{C^0([0,\delta])}.
\end{align}
Finally, we choose $\delta$ small enough (depending on $\varepsilon_1$,
$\|f^0_p(\cdot)-f_{pe}\|_{W^{1,\infty}}$, $\|F\|_{W^{1,\infty}}$) such that
\begin{align}\label{estimatedelta0}
       &\|\mathfrak{F}(\bar\Psi)-\mathfrak{F}(\Psi)\|_{C^0([0,\delta])}\nonumber
         \\
   : =&\max \Big\{ \|\mathfrak{F}_1(\bar\Psi)-\mathfrak{F}_1(\Psi)\|_{C^0([0,\delta])},
                     \|\mathfrak{F}_2(\bar\Psi)-\mathfrak{F}_2(\Psi)\|_{C^0([0,\delta])}  \Big\} \nonumber
            \\
    \leq  & \frac{1}{2}  \|\bar \Psi-\Psi \|_{C^0([0,\delta])},
\end{align}
Banach fixed point  theorem implies the existence of the unique fixed point
$(l(\cdot),f_p(\cdot,1))$ of the mapping $\mathfrak{F}$: $\Psi = \mathfrak{F}(\Psi)$
in $\Omega_{\delta,\vep_1}$.

\noindent \textbf{Step 2. Construction of a solution by characteristic method.}

With the existence of $(l(\cdot),f_p(\cdot,1))$ and $\delta>0$ by \textbf{Step 1},
we can construct a solution to Cauchy problem \eqref{eql-normal}-\eqref{eq-filling-ratio-bound}.
For every $(t,x)$ in $[0,\delta]\times [0,1]$, we still denote by $\xi(s;t,x)$, with $(s, \xi(s;t,x)) \in [0,t]\times [0,1]$,
the characteristic curve passing through the point $(t,x)$; see \e{odexi2} and Fig.\ref{Fig3}.
Since the velocity function $\alpha_p$ is positive, the characteristic $\xi(s;t,x)$ intersects the $x$-axis at point $(0,\beta(t,x))$
with $\beta(t,x)=\xi(0;t,x)$ if  $0\leq  \xi(t;0,0) \leq x\leq 1$;
the characteristic $\xi(s;t,x)$ intersects the $t$-axis at point $(\tau(t,x),0)$ with
$\xi(\tau(t,x);t,x)=0$ if  $0\leq x\leq \xi(t;0,0)$.
Moreover, we have (see \cite[Lemma 3.2 and its proof, Page 90-91]{LiYuBook} for a more general situation)
\begin{align}
\f{\pa \tau(t,x)}{\pa t}
  &=\f{-l(\tau(t,x))}{\zeta N(\tau(t,x))} \f{ \zeta N(t)-xF(l(t),N(t),f_p(t,1) )}{l(t)}
   \,  e^{\int_{\tau(t,x)}^t \f{F(l(\sigma),N(\sigma),f_p(\sigma,1))}{l(\sigma)}\, d\sigma},
   \\
\label{applem12}
\f{\pa \tau(t,x)}{\pa x}&=\f{-l(\tau(t,x))}{\zeta N(\tau(t,x))}
  \, e^{\int_{\tau(t,x)}^t \f{F(l(\sigma),N(\sigma),f_p(\sigma,1))}{l(\sigma)}\, d\sigma},
   \\
\f{\pa \beta(t,x)}{\pa t}
    &= - \f{ \zeta N(t)-xF(l(t),N(t),f_p(t,1) )}{l(t)}
       \, e^{\int_0^t \f{F(l(\sigma),N(\sigma),f_p(\sigma,1))}{l(\sigma)}\, d\sigma},
       \\
          \label{applem22}
\f{\pa \beta(t,x)}{\pa x}
    &= e^{\int_0^t \f{F(l(\sigma),N(\sigma),f_p(\sigma,1))}{l(\sigma)}\, d\sigma}.
\end{align}
We define $f_p$ by
\be\label{sol1}
f_p(t,x)=
\begin{cases}
\displaystyle\frac{F_{in}(\tau(t,x))}{{\color{red} \rho_{o} }V_{eff}N(\tau(t,x))},
            \quad \text{if}\ 0\leq  x \leq  \xi(t;0,0) \leq 1, 0\leq t\leq \delta,
             \\
 f^{0}_p(\beta(t,x)),\quad \text{if}\  0\leq  \xi(t;0,0)\leq  x\leq  1, 0\leq t\leq \delta.
 \end{cases}
\ee
Then it is easy to check that
$(l,f_p)\in W^{1,\infty}(0,\delta) \times W^{1,\infty}((0,\delta)\times (0,1))$
under the compatibility condition \e{compati} and $(l,f_p)$ is indeed a solution
to Cauchy problem \eqref{eql-normal}-\eqref{eq-filling-ratio-bound}.

\noindent \textbf{Step 3. Uniqueness of the solution.}

Assume that Cauchy problem  \eqref{eql-normal}-\eqref{eq-filling-ratio-bound}
has two solutions $(l,f_p),(\bar l,\bar f_p)$ on $[0,\delta] \times [0,1]$.
 It follows that $(l(\cdot),f_p(\cdot,1))=(\bar l(\cdot),\bar f_p(\cdot,1))$
since they are both the fixed point of the mapping $\mathfrak{F}$: $\Psi = \mathfrak{F}(\Psi)$
in $\Omega_{\delta,\vep_1}$. This fact implies that the characteristics $\xi(\cdot;t,x)$ and
 $\bar \xi(\cdot;t,x)$ coincide with each other
and therefore so do the solutions $f_p$ and $\bar f_p$ by characteristic method.

\noindent \textbf{Step 4. A priori estimate on the local solution.}

By definition of $f_p$ and assumption  \eqref{initial-bound}, it is clear that  for all $t\in[0,\delta]$,
\be
\|f_p(t,\cdot)-f_{pe}\|_{L^{\infty}}\leq  \varepsilon. \label{G4}
\ee
Thanks to \e{eql-normal}, \e{equil} ,  \e{G4} and assumption \e{initial-bound},
we get for all $t\in [0,\delta]$ that
\begin{align}
             |\dot{l}(t)|
      =&|F\left(l(t),N(t),f_p(t,1)\right)-F\left(l_e,N_e,f_{pe}\right)|\nonumber
             \\
    \leq & \|F\|_{W^{1,\infty}}  (|l(t)-l_e|+|N(t)-N_e|+|f_p(t,1)-f_{pe}|).  \label{G5}
            \\
     \leq  &   \|F\|_{W^{1,\infty}} |l(t)-l_e|+ 2 \vep \|F\|_{W^{1,\infty}},
   \end{align}
which yields \eqref{estimatel} from \eqref{initial-bound} and Gronwall's inequality.
On the other hand,
\begin{align}\label{partialfpx}
        \Big\|\f{\pa f_p}{\pa x}\Big\|_{L^{\infty}}
   \leq  &\Big\|\f{\pa }{\pa x} \Big(\frac{F_{in}(\tau(t,x))}{{\color{red} \rho_{o} }V_{eff}N(\tau(t,x))}\Big) \Big\|_{L^{\infty}}
                   +\Big\|\f{\pa }{\pa x} f^0_p(\beta(t,x)) \Big\|_{L^{\infty}}\nonumber
           \\
    \leq  &\Big\|\frac{F_{in}(\cdot)}{{\color{red} \rho_{o} }V_{eff}N(\cdot)}-f_{pe}\Big\|_{W^{1,\infty}}
                      \Big\|\f{\pa \tau}{\pa x}\Big\|_{L^{\infty}}
      +\Big\|f^0_p(\cdot)-f_{pe}\Big\|_{W^{1,\infty}} \Big\|\f{\pa \beta}{\pa x}\Big\|_{L^{\infty}}.
\end{align}
Combining \eqref{applem12}, \eqref{applem22}, \eqref{partialfpx} and assumption \eqref{initial-bound},
we obtain \eqref{estimatefp} which concludes
the proof of Lemma \ref{lem-loc}.    \qed


\section{Proof of Theorem \ref{thm-regu} and Theorem \ref{thm-mt}} \label{sec-proof 2}

Before proving Theorem \ref{thm-regu} and Theorem \ref{thm-mt},
let us recall a classical result on Cauchy problem of the following general linear transport equation
\be\label{equation-u}
\begin{cases}
u_t+a(t,x)u_x=b(t,x)u+c(t,x),\quad (t,x)\in Q=(0,T)\times(0,1),\\
u(0,x)=u_0(x),\quad x\in (0,1),\\
u(t,0)=h(t),\quad t\in (0,T),
\end{cases}
\ee
where $a(t,x) > 0$, $a, a_x, b\in L^{\infty}(Q)$ and $c\in L^2(Q)$.

We recall from \cite[Section 2.1]{CoronBook}, the definition of a weak solution to Cauchy problem
\eqref{equation-u}.
\begin{definition}\label{definition-u}
Let $T>0$, $u_0\in L^2(0,1)$, $h\in L^2(0,T)$ be given.
A weak solution of Cauchy problem
\eqref{equation-u} is a function
$u\in C^0([0,T];L^2(0,1))$
such that for every $\tau\in[0,T]$, every test function $\varphi\in C^1([0,T]\times[0,1])$ such that
$ \varphi(t,1)=0,\ \forall t\in[0,T]$,
one has
\begin{align}\label{definition-u1}
&- \int_0^{\tau}\int_0^1 \Big(u [\pa_t\varphi+a \pa_x \varphi+(a_x+b)\varphi ]  +  c\varphi \Big)\, dx\, dt
   + \int_0^1 u(\tau,\cdot) \varphi(\tau,\cdot)  \, dx   \nonumber\\
&- \int_0^1 u_0\varphi(0,\cdot)\, dx
- \int_0^{\tau}ha(\cdot,0)\varphi(\cdot,0)\, dt=0.
\end{align}
\end{definition}
We have the following lemma
\begin{lem}\label{lem-lin}
Let $T>0$, $u_0 \in L^2(0,1)$ and $h \in L^2(0,T)$ be given.
Then, Cauchy problem \eqref{equation-u} has
a unique weak solution $u \in C^0([0,T];L^2(0,1))$ and the following estimate holds:
\begin{align}\label{ape}
\|u\|_{C^0([0,T];L^2(0,1))}\leq  C (\|u_0\|_{L^2(0,1)}+\|h\|_{L^2(0,T)}+\|c\|_{L^2(Q)} ),
\end{align}
where $C=C(T, \|a\|_{L^{\infty}(Q)},\|a_x\|_{L^{\infty}(Q)},\|b\|_{L^{\infty}(Q)})$
is a constant independent of $u_0,h,c$.
\end{lem}

For the proof of Lemma \ref{lem-lin}, one can refer to \cite{LiYuBook} for classical solution or
\cite[Theorem 23.1.2, Page 387]{Horma3} for Cauchy problem on $\mathbb{R}$ without boundary.

\noindent
{\bf Proof of Theorem \ref{thm-regu}.}
By Theorem \ref{thm-well} and  Lemma \ref{lem-lin}, it suffices to prove that  the systems of $f_{p_{xx}}$
satisfies all the assumptions of Lemma \ref{lem-lin}.

Differentiating \eqref{eq-filling-ratio-bound} with respect to $x$ once and twice give us successively  that
\be\label{de-eq-filling-ratio-bound}
\begin{cases}
\partial_t f_{p_x}(t,x)+\alpha_p(t,x) \partial_x f_{p_x}(t,x)=-\alpha_{p_x}(t,x)  f_{p_x}(t,x),\quad (t,x)\in Q,\\
f_{p_x}(0,x)=f^{0}_{p_x}(x),\quad x\in(0,1),\\
f_{p_x}(t,0)=\displaystyle\f{-l(t)}{\zeta N(t)}\cdot \Big(\f{F_{in}(t)}{{\color{red} \rho_{o} }V_{eff}N(t)}\Big)',\quad t\in(0,T),
\end{cases}
\ee
and
\be\label{de-de-eq-filling-ratio-bound}
\begin{cases}
\partial_t f_{p_{xx}}(t,x)+\alpha_p(t,x)  \partial_x f_{p_{xx}}(t,x)
       =-2\alpha_{p_x}(t,x) f_{p_{xx}}(t,x),\quad (t,x)\in Q,
         \\
f_{p_{xx}}(0,x)=f^{0}_{p_{xx}}(x),\quad x\in(0,1),
       \\
f_{p_{xx}}(t,0)=\displaystyle\f{-l(t)}{\zeta N(t)}
       \Big [\f{F(l(t),N(t),f_p(t,1))}{\zeta N(t)} \Big(\f{F_{in}(t)}{{\color{red} \rho_{o} }V_{eff}N(t)}\Big)' \\
       \qquad \qquad \qquad \qquad  \qquad \displaystyle
          -\Big(\f{l(t)}{\zeta N(t)} \Big(\f{F_{in}(t)}{{\color{red} \rho_{o} }V_{eff}N(t)}\Big)'\Big)' \Big ], \quad t\in(0,T),
\end{cases}
\ee
with
     \be  \alpha_p(t,x) =\f{\zeta N(t)-xF(l(t),N(t),f_p(t,1))}{l(t)},
     \quad  \alpha_{p_x}(t,x)=\f{-F(l(t),N(t),f_p(t,1))}{l(t)}.
\ee
From the assumptions that $f_p^0\in H^2(0,1)$, $\displaystyle\f{F_{in}(\cdot)}{{\color{red} \rho_{o} }V_{eff}N(\cdot)}\in H^2(0,T)$
and the compatibility conditions \e{compati} and \e{compati2},
one easily concludes Theorem \ref{thm-regu}   by applying Lemma \ref{lem-lin}
  to Cauchy problem \eqref{de-de-eq-filling-ratio-bound}.              \qed

\noindent
{\bf Proof of Theorem \ref{thm-mt}. }
By Theorem  \ref{thm-well}, we have already $(l,f_p)\in W^{1,\infty}(0,T)\times W^{1,\infty}(Q)$.
Then Theorem \ref{thm-mt} is a direct consequence of Lemma \ref{lem-lin} by solving first $(M_p,T_p)$
in the $PFZ$ and next $(M_f,T_f)$ in the $FFZ$.      \qed

\section{Conclusion}\label{conclusion}

In this paper, we consider the well-posedness of the Cauchy problem for  a physical model
of the extrusion process, which is described by  two systems of conservation laws in
complementary time varying domains.
%
After a suitable change of coordinates in space variables
the original system is transformed into a normalized problem in fixed domain.
Then  we prove the existence and uniqueness of $(l(t),f_p(t,1))$ for
$t$ small by Banach fixed point  theorem. With the known $(l(t),f_p(t,1))$,
we construct a local solution to this Cauchy problem and prove that the local solution is unique.
Using the estimates on the local solution and  induction in time,
we can extend the local solution to the semi-global one. The solution we obtained in Theorem
\ref{thm-well} is called semi-global solution since it exists
on any preassigned time interval $[0,T]$ if the initial and
boundary data has some kind of smallness (depending on $T$).
The $H^2$-regularity of the filling ratio $f_p$ is also proved as preliminaries for
 asymptotic stabilization for the corresponding closed-loop system with feedback controls.

Based on the analytical results obtained in this paper,
we are able to study the controllability and stabilization of this model which is important in applications.
It is interesting, in particular, to study the controllability of boundary profile, i.e., to reach the desired moisture
and temperature at the die under suitable controls.
These interesting control problems will be studied in some forthcoming papers.
We  point out  also that the assumption  of constant fluid viscosity  along
the extruder allows  to decouple the  interface dynamics to  moisture  and temperature of the mixture.
However, for many extruded material  the fluid viscosity may
significantly change  with the chemical composition and temperature
evolutions.   Our future works  also consist in analysis and control
of this extrusion process model  in the case of distributed viscosity.
The study on these problem is really challenging for mathematical analysis but also more useful in applications.




\section{Acknowledgements}
This work has been done in the context of the PhD  of Mamadou Diagne at LAGEP-UCBL (Laboratoire d'automatique et du G\'{e}nie des Proc\'{e}d\'{e}s de l'Universit\'{e} Claude Bernard Lyon 1) under the supervision of Professor Bernhard Maschke.

The authors would like to thank Professor Jean-Michel Coron and Professor Bernhard Maschke for their helpful comments and constant support.
The authors are thankful to  the support  of the ERC advanced grant 266907 (CPDENL) and the hospitality of the Laboratoire Jacques-Louis Lions of Universit\'{e} Pierre et Marie Curie.

Peipei Shang was partially supported by the National Science Foundation of China (No. 11301387). Mamadou Diagne  has been supported by a doctoral Grant of the French Ministry of Higher Education and Research and in the context of the French National Research Agency sponsored project ANR-11-BS03-0002 HAMECMOPSYS.
Zhiqiang Wang was partially supported by the Natural Science Foundation of Shanghai
(No. 11ZR1402500) and by the National Science Foundation of China (No. 11271082).

\bibliographystyle{plain}

\bibliography{extruder}

\end{document}